# Existence and concentration of solutions for a class of biharmonic equations

Marcos T.O. Pimenta [a,1], Sérgio H.M. Soares [b,*,2]

[a] Departamento de Matemática, Universidade Estadual de Londrina, 86051-990, Londrina, PR, Brazil
[b] Departamento de Matemática, Instituto de Ciências Matemáticas e de Computação, Universidade de São Paulo, 13560-970, São Carlos, SP, Brazil

ABSTRACT

Some superlinear fourth order elliptic equations are considered. A family of solutions is proved to exist and to concentrate at a point in the limit. The proof relies on variational methods and makes use of a weak version of the Ambrosetti–Rabinowitz condition. The existence and concentration of solutions are related to a suitable truncated equation.

## 1. Introduction

In the last years, many authors have been studied several questions about the following Schrödinger elliptic equation

$$\begin{cases} -\epsilon^2 \Delta u + V(x)u = f(x,u) & \text{in } \Omega, \\ u \in H^1(\Omega) \end{cases} \tag{1.1}$$

with Neumann or Dirichlet boundary conditions, where $\Omega$ is a domain in $\mathbb{R}^N$. Motivated by Floer and Weinstein [12], Rabinowitz in [16] uses a mountain-pass type argument to find ground-state solutions to (1.1) for $\epsilon > 0$ sufficiently small, when $N \geqslant 3$, $\Omega = \mathbb{R}^N$, $f$ is a subcritical and superlinear nonlinearity function and the potential $V$ is nonnegative and assumed to satisfy the condition

$$0 < V_0 = \inf_{x \in \mathbb{R}^N} V(x) < \liminf_{|x| \to \infty} V(x). \tag{1.2}$$

In [20], Wang proves that the mountain-pass solutions found in [16] concentrate around a global minimum of $V$ as $\epsilon \to 0$. In [2,3], Alves and Figueiredo consider the problem (1.1) for the p-Laplace operator obtaining existence, multiplicity and concentration of positive solutions. In the celebrated paper [10], del Pino and Felmer obtained existence and concentration of solutions for the problem (1.1), where $N \geqslant 3$, $f$ is a subcritical and superlinear nonlinearity function and the potential $V$ is nonnegative and it is assumed to satisfy the following condition

$$\inf_{x \in \Lambda} V(x) < \inf_{x \in \partial\Lambda} V(x),$$

where $\Lambda$ is a bounded domain compactly contained in $\Omega$. They developed a penalization-type method in order to overcome the lack of compactness and used the Mountain Pass Theorem to get existence and concentration of solutions. These

* Corresponding author. Fax: +55 16 3373 9650.
*E-mail addresses:* mtopimenta@uel.br (M.T.O. Pimenta), monari@icmc.usp.br (S.H.M. Soares).

[1] Research supported by CNPq, Brazil.
[2] The author was partially supported by CNPq, Brazil.

0022-247X/$ – see front matter

arguments have inspired many authors in the last years, among them we could cite [4] and [11], where the authors have obtained multiplicity and concentration of nodal and positive solutions, respectively, to an equation related to (1.1). In [8], Alves and Soares obtain existence and concentration of nodal solutions of (1.1) for the case where $N = 2$ and the function $f$ has critical exponential growth.

Although there are many works dealing with problem (1.1) and with related p-Laplacian ones, just a few works can be found dealing with biharmonic or even polyharmonic Schrödinger equations. Among them we could cite [5] and [6], where the authors have obtained nontrivial solutions to semilinear biharmonic problems with critical nonlinearities and also [18], where the authors obtained infinitely many solutions for a polyharmonic Schrödinger equation with non-homogeneous boundary data on unbounded domains.

Motivated by the results just described, a natural question is whether the same phenomenon of concentration occurs for the following class of fourth order elliptic equations

$$\begin{cases} \epsilon^4 \Delta^2 u + V(x)u = f(u) & \text{in } \mathbb{R}^N, \\ u \in H^2(\mathbb{R}^N), \end{cases} \tag{1.3}$$

where $\Delta^2$ is the biharmonic operator, $\epsilon > 0$ and $N \geqslant 5$. The functions $f : \mathbb{R} \to \mathbb{R}$ and $V : \mathbb{R}^N \to \mathbb{R}$ satisfy the following assumptions:

$(V_1)$ $V \in C^0(\mathbb{R}^N) \cap L^\infty(\mathbb{R}^N)$.

$(V_2)$ There exist a bounded domain $\Omega \subset \mathbb{R}^N$ and $x_0 \in \Omega$, such that

$$0 < V(x_0) = V_0 = \inf_{\mathbb{R}^N} V < \inf_{\partial\Omega} V.$$

$(f_1)$ $f \in C^1(\mathbb{R})$.

$(f_2)$ $f(0) = f'(0) = 0$.

$(f_3)$ There exist constants $c_1, c_2 > 0$ and $p \in (1, 2_* - 1)$, such that

$$\left|f(s)\right| \leqslant c_1|s| + c_2|s|^p, \quad \forall s \in \mathbb{R},$$

where $2_* = 2N/(N-4)$.

$(f_4)$ $\lim_{|s|\to\infty} \frac{F(s)}{s^2} = +\infty$, where $F(s) = \int_0^s f(t)\,dt$.

$(f_5)$ The function $f(s)/s$ is increasing for $s > 0$ and decreasing for $s < 0$.

Our main result is the following:

**Theorem 1.1.** *Assume that conditions $(V_1)$, $(V_2)$ and $(f_1)$–$(f_5)$ hold. Then for each sequence $\epsilon_n \to 0$, there exists a subsequence, still denoted by $\{\epsilon_n\}$, such that, for all $n \in \mathbb{N}$, there exists a nontrivial weak solution $u_n$ of (1.3) (with $\epsilon = \epsilon_n$). Moreover, if $x_n$ is the maximum point of $|u_n|$, then $x_n \in \Omega$ and*

$$\lim_{n\to\infty} V(x_n) = \inf_{\mathbb{R}^N} V.$$

Although the principal arguments used here can be found in [10], the proofs have to be deeply modified because of some natural difficulties that the study of the biharmonic operator gives rise. For instance, in [10], in order to prove that the solutions of the penalized problem in fact are solutions of the original one, the authors use an argument that relies on the strong maximum principle to the Laplace operator and also on the fact that $u^+ = \max\{u, 0\}$ and $u^- = \min\{u, 0\}$ belong to $H^1(\Omega)$ for every $u$ in $H^1(\Omega)$. In [7], to prove the same to a quasilinear problem, the authors combine a comparison principle together with Moser's iteration technique. However, all these arguments have severe limitations to deal with the biharmonic equation, because of the lack of a general form of the maximum principle to the biharmonic operator and the impossibility of splitting $u = u^+ + u^-$ in $H^2(\Omega)$. To overcome these problems, our argument relies on proving that rescalings of solutions of the penalized problem exhibit a uniform decay in infinity. To prove this we use some compactness results in Nehari manifolds found in [7] together with *a priori* $L^p$ estimates found in [1] and $L^\infty$ estimates proved by Ramos in [15].

It is worth to point out that we provide our results assuming a weaker version of the famous Ambrosetti–Rabinowitz condition (see $(f_4)$). The use of this weaker condition brings some difficulty to prove that the $(PS)$ sequences are bounded, which required some arguments found in [13] (see also [17]). Moreover, this weak condition represents a difficulty to prove that the Nehari manifold is homeomorphic to the unitary sphere in $H^2(\mathbb{R}^N)$. This last problem can be dropped out using similar arguments of Szulkin and Weth in [19].

The article is organized in the following way: In the second section, we use the argument given by [10] and [13] to modify the function $f$ to get the Palais–Smale condition for the functional associated with the respective modified equation. The existence and concentration of solutions to the modified problem are established. Finally, in the third section we prove that these solutions have a kind of uniform decay at infinity.

## 2. Preliminary results

We start by observing that the following problem

$$\begin{cases} \Delta^2 v + V(\epsilon x)v = f(v) & \text{in } \mathbb{R}^N, \\ v \in H^2(\mathbb{R}^N) \end{cases} \tag{2.1}$$

is equivalent to (1.3). In fact, the solutions $v_\epsilon$ of (2.1) and $u_\epsilon$ of (1.3) are related by

$$v_\epsilon(x) = u_\epsilon(\epsilon x).$$

Let $V_0$ be as in $(V_1)$ and let us choose $k > 0$ such that $k > 2V_0$. Let $a > 0$ be a number such that $\max\{\frac{f(a)}{a}, \frac{f(-a)}{-a}\} \leqslant \frac{V_0}{k}$. Set

$$\tilde{f}(s) = \begin{cases} \frac{-f(-a)}{a}s & \text{if } s < -a, \\ f(s) & \text{if } |s| \leqslant a, \\ \frac{f(a)}{a}s & \text{if } s > a. \end{cases}$$

By the continuity of $V$, there exists a non-empty open set $\Omega' \subset \Omega$, such that:

$$\inf_{\Omega \backslash \Omega'} V > \inf_{\mathbb{R}^N} V \quad \text{and} \quad \min_{\partial\Omega'} V > \inf_{\mathbb{R}^N} V.$$

Let $\chi \in C^\infty(\mathbb{R}^N)$, $0 \leqslant \chi \leqslant 1$, be such that

$$\chi(x) = \begin{cases} 1 & \text{if } x \in \Omega', \\ 0 & \text{if } x \in \mathbb{R}^N \backslash \Omega, \end{cases}$$

and define

$$g(x,s) = \chi(x)f(s) + \big(1 - \chi(x)\big)\tilde{f}(s). \tag{2.2}$$

By $(f_1)$–$(f_5)$, $g$ satisfies

$(g_1)$ $g(x,s) = o(|s|)$ as $s \to 0$.
$(g_2)$ There exist $c_1, c_2 > 0$ and $p \in (1, 2_* - 1)$, such that $|g(x,s)| \leqslant c_1|s| + c_2|s|^p$, for all $s \in \mathbb{R}$ and $x \in \mathbb{R}^N$.
$(g_3)$

$$2G(x,s) \leqslant g(x,s)s, \quad \text{for all } x \in \mathbb{R}^N \text{ and } s \in \mathbb{R},$$
$$g(x,s)s \leqslant \frac{1}{k}V(x)s^2, \quad \text{for all } x \notin \Omega \text{ and } s \in \mathbb{R}.$$

$(g_4)$ $g(x,s)/s$ is nondecreasing for $s > 0$ and nonincreasing for $s < 0$, where $x \in \mathbb{R}^N$.

The problem we now consider is the following:

$$\begin{cases} \Delta^2 v + V(\epsilon x)v = g(\epsilon x, v) & \text{in } \mathbb{R}^N, \\ v \in H^2(\mathbb{R}^N). \end{cases} \tag{2.3}$$

Let $E_\epsilon = (H^2(\mathbb{R}^N), \langle\cdot,\cdot\rangle_\epsilon)$ be the Hilbert space endowed with the inner product

$$\langle u, v\rangle_\epsilon = \int_{\mathbb{R}^N} \big(\Delta u \Delta v + V(\epsilon x)uv\big)\,dx.$$

Denote by $\|\cdot\|_\epsilon$ the norm associated with this inner product. We consider the functional $I_\epsilon$ defined on $E_\epsilon$ by

$$I_\epsilon(u) = \frac{1}{2}\int_{\mathbb{R}^N} \big(|\Delta u|^2 + V(\epsilon x)u^2\big)\,dx - \int_{\mathbb{R}^N} G(\epsilon x, u)\,dx,$$

where $G(x,s) = \int_0^s g(x,t)\,dt$. The functional $I_\epsilon \in C^1(E_\epsilon, \mathbb{R}^N)$ and

$$I'_\epsilon(u)v = \int_{\mathbb{R}^N} \big(\Delta u \Delta v + V(\epsilon x)uv\big)\,dx - \int_{\mathbb{R}^N} g(x,u)v\,dx,$$

for all $u, v \in E_\epsilon$. Hence, critical points of $I_\epsilon$ are weak solutions of the Euler–Lagrange equation (2.3).

Our first lemma provides conditions under which $I_\epsilon$ satisfies the hypotheses of the Mountain Pass Theorem.

**Lemma 2.1.** *Assume that conditions $(g_1)$–$(g_3)$ and $(V_1)$ hold. Then, for each $\epsilon > 0$, there exist positive constants $\rho, \beta$ and $\phi \in E_\epsilon$ with $\|\phi\|_\epsilon > \rho$, such that*

1. $I_\epsilon(u) \geqslant \beta$ *for all* $\|u\|_\epsilon = \rho$.
2. $I_\epsilon(\phi) < 0$.

**Proof.** Using $(g_1)$ and $(g_2)$ and the Sobolev imbedding, we can prove that for all $\eta > 0$, there exists a constant $C(\eta) > 0$ such that

$$\int_{\mathbb{R}^N} \left|G(\epsilon x, u)\right| dx \leqslant \eta \|u\|_\epsilon^2 + C(\eta)\|u\|_\epsilon^{p+1}.$$

Hence, by choosing $\eta \in (0, 1/2)$, there exists a small $\rho > 0$ such that

$$I_\epsilon(u) \geqslant \beta > 0, \quad \text{for all } \|u\|_\epsilon = \rho,$$

where $\rho = [(1/2 - \eta) - C(\eta)r^{p-1}]r^2$. This establishes 1.

In order to prove 2, fix $\varphi \in C_0^\infty(\Omega'_\epsilon)$ with $\varphi > 0$, where $\Omega'_\epsilon = \epsilon^{-1}\Omega'$. By $(f_4)$, for every $M \geqslant \|\varphi\|_\epsilon^2/2\|\varphi\|_{L^2}^2$, there exists a constant $c_0 > 0$ such that

$$F(s) \geqslant M|s|^2 - c_0, \quad \text{for all } s \in \mathbb{R}.$$

Then,

$$\begin{aligned} I_\epsilon(t\varphi) &= \frac{t^2}{2}\|\varphi\|_\epsilon^2 - \int_{\mathbb{R}^N} F(t\varphi)\, dx \\ &\leqslant \frac{t^2}{2}\|\varphi\|_\epsilon^2 - t^2 M \int_{\mathbb{R}^N} |\varphi|^2\, dx + c_0\left|supp(\varphi)\right| \\ &= t^2\left(\frac{\|\varphi\|_\epsilon^2}{2} - M\int_{\mathbb{R}^N} |\varphi|^2\, dx\right) + c_0\left|supp(\varphi)\right|. \end{aligned}$$

Therefore, $I_\epsilon(t\varphi) \to -\infty$ as $t \to +\infty$ and the proof is complete. □

**Lemma 2.2.** *Assume that conditions $(g_1)$–$(g_3)$ and $(V_1)$ hold. Then, the functional $I_\epsilon$ satisfies the Palais–Smale condition, that is, if $\{u_n\}$ is a sequence in $E_\epsilon$ such that $\{I_\epsilon(u_n)\}$ is bounded and $I'_\epsilon(u_n) \to 0$, then $\{u_n\}$ contains a strongly convergent subsequence in $E_\epsilon$.*

**Proof.** We start by claiming that $\{u_n\}$ is a bounded sequence in $E_\epsilon$. Let us suppose that $I_\epsilon(u_n) \to d$ and $I'(u_n)u_n \to 0$. Suppose by contradiction that $\|u_n\|_\epsilon \to +\infty$. Let $w_n = u_n/\|u_n\|_\epsilon$. Analogously to the proof of Lemma 2.5 below, one can prove that $w_n \rightharpoonup 0$ in $E_\epsilon$. Let $t_n \in [0, 1]$ be such that

$$I_\epsilon(t_n w_n) = \max_{t\in[0,1]} I_\epsilon(t w_n).$$

Letting $n_0 \in \mathbb{N}$ be such that $R/\|u_n\|_\epsilon \leqslant 1$ for all $n \geqslant n_0$, using the definition of $t_n$ and the condition $(g_3)$, it follows that

$$\begin{aligned} I_\epsilon(t_n u_n) &\geqslant I_\epsilon\left(\frac{R}{\|u_n\|_\epsilon} u_n\right) \\ &= \frac{R^2}{2}\|w_n\|_\epsilon^2 - \int_{\Omega_\epsilon} G(\epsilon x, R w_n)\, dx - \int_{\mathbb{R}^N\setminus\Omega_\epsilon} G(\epsilon x, R w_n)\, dx \\ &\geqslant \frac{R^2}{2}\|w_n\|_\epsilon^2 - \int_{\Omega_\epsilon} G(\epsilon x, R w_n)\, dx - \frac{R^2}{2}\int_{\mathbb{R}^N} \frac{V(\epsilon x)}{k} w_n^2\, dx \\ &\geqslant \frac{R^2(k-1)}{2k} - \int_{\Omega_\epsilon} G(\epsilon x, R w_n)\, dx. \end{aligned}$$

Since $w_n \rightharpoonup 0$ in $E_\epsilon$ we have that $w_n \to 0$ in $L^r_{loc}(\mathbb{R}^N)$ for all $2 < r < 2_*$. By $(g_2)$,

$$\liminf_{n\to+\infty} I(t_n u_n) \geqslant \frac{R^2(k-1)}{2k}$$

and since $R > 0$ is arbitrary we get

$$\lim_{n\to+\infty} I(t_n u_n) = +\infty.$$

This together with the fact that $I_\epsilon(0) = 0$ and $I_\epsilon(u_n) \to d$ implies that $t_n \in (0,1)$. Therefore

$$I'_\epsilon(t_n u_n)u_n = 0.$$

Since by $(g_4)$, $H(x,s) = g(x,s)s - 2G(x,s)$ is nondecreasing in $s > 0$ and nonincreasing in $s < 0$, we have

$$2I_\epsilon(t_n u_n) = \int_{\mathbb{R}^N} \big(g(\epsilon x, t_n u_n)t_n u_n - 2G(\epsilon x, t_n u_n)\big)\,dx$$
$$\leqslant \int_{\mathbb{R}^N} \big(g(\epsilon x, u_n)u_n - 2G(\epsilon x, u_n)\big)\,dx,$$

which implies that

$$\lim_{n\to+\infty} \int_{\mathbb{R}^N} \big(g(\epsilon x, u_n)u_n - 2G(\epsilon x, u_n)\big)\,dx = +\infty. \tag{2.4}$$

On the other hand,

$$\int_{\mathbb{R}^N} \big(g(\epsilon x, u_n)u_n - 2G(\epsilon x, u_n)\big)\,dx = 2I_\epsilon(u_n) - I'_\epsilon(u_n)u_n = 2d + o_n(1),$$

which contradicts (2.4) and proves the claim. Then, we can assume that

$$\begin{aligned} &u_n \rightharpoonup u \quad \text{in } E_\epsilon,\\ &u_n \to u \quad \text{in } L^q_{loc}\big(\mathbb{R}^N\big), \text{ for all } 1 \leqslant q < 2_*,\\ &u_n \to u \quad \text{a.e. in } \mathbb{R}^N \end{aligned} \tag{2.5}$$

as $n \to \infty$. By Lebegue's convergence theorem, it is a simple matter to verify that $u$ is a weak solution of (2.3). We now take advantage of the Hilbertian structure of $E_\epsilon$ to prove that $u_n \to u$, as $n \to \infty$, by proving that $\|u_n\|_\epsilon \to \|u\|_\epsilon$ as $n \to \infty$. As we will see, this follows from the following claim.

**Claim.** *Given $\delta > 0$, there exists $R = R(\delta) > 0$ such that*

$$\limsup_{n\to\infty} \int_{\mathbb{R}^N\setminus B_R} \big(|\Delta u_n|^2 + V(\epsilon x)u_n^2\big)\,dx < \delta.$$

Effectively, for each $R > 0$ let $\eta_R \in C^\infty(\mathbb{R}^N)$ be a cut-off function such that $0 \leqslant \eta_R \leqslant 1$, $\eta_R = 0$ in $B_{R/2}(0)$, $\eta_R = 1$ in $B_R(0)^c$, $|\nabla\eta_R| \leqslant C/R$ and $|\Delta\eta_R| \leqslant C/R^2$. From $(g_3)$ and the Hölder inequality, for $R > 0$ such that $\Omega_\epsilon \subset B_{R/2}(0)$, we have

$$\left(1 - \frac{1}{k}\right) \int_{B_R(0)^c} \big(|\Delta u_n|^2 + V(\epsilon x)u_n^2\big)\,dx \leqslant I'_\epsilon(u_n)(\eta_R u_n) + \frac{C}{R},$$

and the claim follows taking the supremum limit.

Combining the claim with the Sobolev imbedding theorem and the integrability of $x \mapsto g(\epsilon x, u(x))u(x)$, we have that given $\delta > 0$ there exists $R_\delta > 0$ such that

$$\limsup_{n\to\infty} \int_{B_R(0)^c} g(\epsilon x, u_n)u_n\,dx < \frac{\delta}{2},$$

and

$$\int_{B_R(0)^c} g(\epsilon x, u)u\,dx < \frac{\delta}{2},$$

for all $R > R_\delta$. By (2.5) and $(g_2)$, we get

$$\left|\int_{\mathbb{R}^N} \big(g(\epsilon x, u_n)u_n - g(\epsilon x, u)u\big)\,dx\right| < \delta,$$

provided $n$ is sufficiently large. Therefore,

$$\|u_n\|_\epsilon - \|u\|_\epsilon = \int_{\mathbb{R}^N} \big(g(\epsilon x, u_n)u_n - g(\epsilon x, u)u\big)\,dx + o_n(1) < \delta + o_n(1)$$

for all $\delta > 0$. The proof of Lemma 2.2 is complete. □

By the Moutain Pass Theorem [9], for any $\epsilon > 0$ there exists $v_\epsilon \in E_\epsilon$ a weak solution of (2.3) such that $I_\epsilon(v_\epsilon) = c_\epsilon$, where

$$c_\epsilon = \inf_{\gamma \in \Gamma_\epsilon} \max_{t\in[0,1]} I_\epsilon\big(\gamma(t)\big)$$

and $\Gamma_\epsilon = \{\gamma \in C^0([0,1], E_\epsilon);\ \gamma(0) = 0 \text{ and } I_\epsilon(\gamma(1)) < 0\}$. From $(g_4)$, the minimax level $c_\epsilon$ can be characterized as (see [16])

$$c_\epsilon = \inf_{u\in E_\epsilon\setminus\{0\}} \max_{t\geqslant 0} I_\epsilon(tu) = \inf_{\mathcal{N}_\epsilon} I_\epsilon,$$

where $\mathcal{N}_\epsilon$ is defined by

$$\mathcal{N}_\epsilon = \big\{u \in E_\epsilon\setminus\{0\};\ I'_\epsilon(u)u = 0\big\}.$$

We note that unlike in [16], there is an additional difficult to prove that $\mathcal{N}_\epsilon$ is homeomorphic to the unitary sphere in $E_\epsilon$ when $f$ does not satisfy the Ambrosetti–Rabinowitz condition. However, the proof of this fact under the condition $(f_4)$ proceeds along the same lines as in [19].

We now consider a sequence $\{\epsilon_n\}$ with $\epsilon_n \to 0$ as $n \to \infty$. We claim that there exists a subsequence, still denoted by $\{\epsilon_n\}$, such that $v_n := v_{\epsilon_n}$ is a solution of (2.1). The proof will be carried out by a series of lemmas. The first one states the existence of a ground-state solution to the limit problem.

**Lemma 2.3.** *Suppose that $f$ satisfies $(f_1)$–$(f_5)$. Then, there exists a ground-state solution to the following problem*

$$\Delta^2 w + V_0 w = f(w) \quad \text{in } \mathbb{R}^N, \tag{2.6}$$

*at the level*

$$c_0 = \inf_{\gamma\in\Gamma_0} \sup_{0\leqslant t\leqslant 1} I_0\big(\gamma(t)\big),$$

*where $I_0$ is the energy functional associated to* (2.6) *and*

$$\Gamma_0 = \big\{\gamma \in C^0\big([0,1], H^2\big(\mathbb{R}^N\big)\big);\ \gamma(0) = 0 \text{ and } I_0\big(\gamma(1)\big) < 0\big\}.$$

**Proof.** Let $E_0$ be the space $H^2(\mathbb{R}^N)$ endowed with the following norm $\|u\|_0^2 = \int_{\mathbb{R}^N}(|\Delta u|^2 + V_0 u^2)\,dx$. In the same way as that in Lemma 2.1 we can prove that $I_0$ satisfies the geometric conditions of Mountain Pass Theorem. Then there exists $\{w_n\} \subset E_0$, such that

$$I_0(w_n) \to c_0 \quad \text{and} \quad I'_0(w_n) \to 0, \quad \text{as } n \to \infty.$$

Using the same arguments as that in Lemma 2.2, one can prove that $\{w_n\}$ is a bounded sequence in $E_0$. By Lions's lemma, it follows that there exist $\{y_n\} \subset \mathbb{R}^N$, $R, \beta > 0$ such that

$$\liminf_{n\to\infty} \int_{B_R(y_n)} w_n^2\,dx > \beta.$$

Defining $u_n(x) = w_n(x + y_n)$, note that $(u_n)$ is bounded in $E_0$. Since $I_0$ is invariant by translations, it follows that $I_0(u_n) \to c_0$, $I'_0(u_n) \to 0$ and

$$\liminf_{n\to\infty} \int_{B_R(0)} u_n^2\,dx > \beta.$$

Hence $u_n$ converges weakly to $u \in E_0\setminus\{0\}$ and strongly in $L^r_{loc}(\mathbb{R}^N)$ for $2 < r < 2_*$, where $u$ is a weak nontrivial solution of (2.6). What is left is to show that $I_0(u) = c_0$. Trivially we have that $c_0 \leqslant I_0(u)$. By the other inequality, let $\rho > 0$. Note that

$$I_0(u_n) - \frac{1}{2}I_0'(u_n)u_n = \int_{\mathbb{R}^N} \left(\frac{1}{2}f(u_n)u_n - F(u_n)\right)dx$$
$$\geqslant \int_{B_\rho(0)} \left(\frac{1}{2}f(u_n)u_n - F(u_n)\right)dx$$

and then, by Fatou's lemma it follows that

$$c_0 \geqslant \int_{B_\rho(0)} \left(\frac{1}{2}f(u)u - F(u)\right)dx.$$

Since the last inequality holds for all $\rho > 0$, it follows by the Lebesgue Dominated Convergence Theorem that

$$c_0 \geqslant \int_{\mathbb{R}^N} \left(\frac{1}{2}f(u)u - F(u)\right)dx = I_0(u)$$

and this implies that $I_0(u) = c_0$. $\square$

**Lemma 2.4.** *Suppose that g satisfies* $(g_1)$–$(g_4)$ *and V satisfies* $(V_1)$–$(V_2)$. *Then*

$$\limsup_{n\to\infty} c_{\epsilon_n} \leqslant c_0.$$

**Proof.** We assume without loss of generality that $x_0 = 0$, for $x_0$ given by condition $(V_2)$. Let $w$ be a solution of (2.6) such that $I_0(w) = c_0$. Let $\psi \in C^\infty(\mathbb{R}^N)$ be a cut-off function such that $\psi \equiv 1$ in $B_\rho(0)$ and $\psi \equiv 0$ in $B_{2\rho}^c(0)$, where $B_{2\rho}(0) \subset \Omega'$. Set $w_n(x) = \psi(\epsilon_n x)w(x)$ and note that $supp(w_n) \subset \Omega'_{\epsilon_n}$, $w_n \to w$ in $H^2(\mathbb{R}^N)$ and in $L^r(\mathbb{R}^N)$ where $2 < r < 2_*$. Let $\varphi_{\epsilon_n}(w_n) > 0$ be such that $\varphi_{\epsilon_n}(w_n)w_n \in \mathcal{N}_{\epsilon_n}$. Suppose that $\varphi_{\epsilon_n}(w_n) \to 1$ as $n \to \infty$. Then

$$c_{\epsilon_n} \leqslant I_{\epsilon_n}\big(\varphi_{\epsilon_n}(w_n)w_n\big)$$
$$= I_0\big(\varphi_{\epsilon_n}(w_n)w_n\big) + \frac{1}{2}\int_{\mathbb{R}^N} \big(V(\epsilon_n x) - V(0)\big)\big(\varphi_{\epsilon_n}(w_n)w_n\big)^2 dx,$$

and the result follows by the Lebesgue Dominated Convergence Theorem. It remains to prove that $\varphi_\epsilon(w_\epsilon) \to 1$ as $\epsilon \to 0$. Since $I'_\epsilon(\varphi_\epsilon(w_\epsilon)w_\epsilon)w_\epsilon = 0$, it follows that

$$\varphi_\epsilon(w_\epsilon)\int_{\mathbb{R}^N} \big(|\Delta w_\epsilon|^2 + V(\epsilon x)w_\epsilon^2\big)dx = \int_{\mathbb{R}^N} f\big(\varphi_\epsilon(w_\epsilon)w_\epsilon\big)w_\epsilon\,dx.$$

We claim that $\{\varphi_\epsilon(w_\epsilon)\}$ is bounded. In fact, on the contrary, there exists $\epsilon_n \to 0$ such that $\varphi_{\epsilon_n}(w_{\epsilon_n}) \to +\infty$. Let $\Sigma \subset \mathbb{R}^N$ be such that $|\Sigma| > 0$ and $w(x) \neq 0$ for all $x \in \Sigma$. Since $H(s) = f(s)s - 2F(s)$ is positive, it holds for all $n \in \mathbb{N}$ that

$$\|w_{\epsilon_n}\|_{\epsilon_n}^2 = \int_{\mathbb{R}^N} \frac{f(\varphi_{\epsilon_n}(w_{\epsilon_n})w_{\epsilon_n})\varphi_{\epsilon_n}(w_{\epsilon_n})w_{\epsilon_n}}{\varphi_{\epsilon_n}(w_{\epsilon_n})^2}dx$$
$$\geqslant \int_\Sigma \frac{2F(\varphi_{\epsilon_n}(w_{\epsilon_n})w_{\epsilon_n})}{\varphi_{\epsilon_n}(w_{\epsilon_n})^2}dx$$
$$= \int_{\Sigma\setminus w_{\epsilon_n}^{-1}(0)^c} \frac{2F(\varphi_{\epsilon_n}(w_{\epsilon_n})w_{\epsilon_n})}{(\varphi_{\epsilon_n}(w_{\epsilon_n})w_{\epsilon_n})^2}w_{\epsilon_n}^2\,dx.$$

On the other hand, by $(f_4)$ and Fatou's lemma it follows that

$$\liminf_{n\to\infty} \int_{\Sigma\setminus w_{\epsilon_n}^{-1}(0)^c} \frac{2F(\varphi_{\epsilon_n}(w_{\epsilon_n})w_{\epsilon_n})}{(\varphi_{\epsilon_n}(w_{\epsilon_n})w_{\epsilon_n})^2}w_{\epsilon_n}^2\,dx = +\infty,$$

which implies that

$$\|w_{\epsilon_n}\|_{\epsilon_n}^2 \to +\infty, \quad \text{as } n \to \infty,$$

which contradicts the fact that $w_{\epsilon_n} \to w$ as $n \to \infty$. We can now verify that $\varphi_\epsilon(w_\epsilon) \nrightarrow 0$ as $\epsilon \to 0$. In fact, on the contrary there exists $\epsilon_n \to 0$ such that $\varphi_{\epsilon_n}(w_{\epsilon_n}) \to 0$ as $n \to \infty$. By $(f_2)$–$(f_3)$ and the Sobolev imbedding one can prove that

$$\lim_{n\to\infty} \int_{\mathbb{R}^N} \frac{f(\varphi_{\epsilon_n}(w_{\epsilon_n})w_{\epsilon_n})w_{\epsilon_n}^2}{\varphi_{\epsilon_n}(w_{\epsilon_n})w_{\epsilon_n}}\,dx = 0. \tag{2.7}$$

On the other hand,

$$\|w_{\epsilon_n}\|_{\epsilon_n}^2 = \int_{\mathbb{R}^N} \frac{f(\varphi_{\epsilon_n}(w_{\epsilon_n})w_{\epsilon_n})w_{\epsilon_n}^2}{\varphi_{\epsilon_n}(w_{\epsilon_n})w_{\epsilon_n}}\,dx. \tag{2.8}$$

Hence by (2.7) and (2.8), one can see that $\|w_{\epsilon_n}\|_{\epsilon_n} \to 0$, which contradicts the fact that $w_{\epsilon_n} \to w$ and $I_{V_0}(w) = c_{V_0} > 0$. Then there exist $\alpha, \beta > 0$ such that

$$\alpha \leqslant \varphi_\epsilon(w_\epsilon) \leqslant \beta.$$

Using that $w_{\epsilon_n} \to w$ in $H^2(\mathbb{R}^N)$ and $w$ is a solution of (2.6), it follows by $(f_5)$ that $\varphi_\epsilon(w_\epsilon) \to 1$. $\square$

The proof of the next result is based upon ideas found in [13].

**Lemma 2.5.** *$\{v_n\}$ is a bounded sequence in $H^2(\mathbb{R}^N)$.*

**Proof.** By Lemma 2.4, we have that $\{I_{\epsilon_n}(v_n)\}$ is bounded and $I'_{\epsilon_n}(v_n) = 0$ for all $n \in \mathbb{N}$. Let us suppose that $\{v_n\} \subset H^2(\mathbb{R}^N)$ is such that $I_{\epsilon_n}(v_n) \to d \leqslant c_0$. Assume by contradiction that $\{v_n\}$ is such that $\|v_n\|_{H^2(\mathbb{R}^N)} \to +\infty$ along a subsequence. Let $w_n = v_n/\|v_n\|_{\epsilon_n}$. Note that $w$ satisfies the following problem

$$\Delta^2 w_n + V(\epsilon_n x)w_n = \chi(\epsilon_n x)\frac{f(v_n)}{v_n}w_n + \big(1-\chi(\epsilon_n x)\big)\frac{\tilde{f}(v_n)}{v_n}w_n \quad \text{in } \mathbb{R}^N. \tag{2.9}$$

We claim that one of the following conditions holds:

i) $$\limsup_{n\to\infty} \sup_{y\in\mathbb{R}^N} \int_{B_1(y)} \big|\chi(\epsilon_n x)w_n(x)\big|^2\,dx > 0;$$

ii) $$\lim_{n\to\infty} \sup_{y\in\mathbb{R}^N} \int_{B_1(y)} \big|\chi(\epsilon_n x)w_n(x)\big|^2\,dx = 0.$$

However, we will show that both conditions in fact do not occur and this will give us a contradiction. Suppose that i) holds. Then there exists $\{y_n\} \subset \mathbb{R}^N$ such that

$$\limsup_{n\to\infty} \int_{B_1(y_n)} \big|\chi(\epsilon_n x)w_n(x)\big|^2\,dx > 0. \tag{2.10}$$

Then note that $B_1(y_n) \cap \Omega_{\epsilon_n} \neq \emptyset$ and we can suppose that $\epsilon_n y_n \to x_0 \in \overline{\Omega}$. Let $\bar{w}_n(x) := w_n(x + y_n)$ and note that $\{\bar{w}_n\}$ is bounded in $H^2(\mathbb{R}^N)$. Then there exists $w_0 \in H^2(\mathbb{R}^N)$ such that $\bar{w}_n \rightharpoonup w_0$ in $H^2(\mathbb{R}^N)$. Hence,

$$\chi(\epsilon_n \cdot + \epsilon_n y_n)\bar{w}_n(\cdot) \rightharpoonup \chi(x_0)w_0, \quad \text{in } H^2\big(\mathbb{R}^N\big).$$

By (2.10),

$$\int_{B_1(0)} \big|\chi(x_0)w_0\big|^2\,dx > 0, \tag{2.11}$$

which implies that $\chi(x_0) \neq 0$ and there exists $\Gamma \subset B_1(0)$, $|\Gamma| > 0$, such that

$$w_0(x) \neq 0 \quad \text{for all } x \in \Gamma. \tag{2.12}$$

By multiplying (2.9) by $w_n$ and integrating by parts, we have that

$$1 = \int_{\mathbb{R}^N} \left(\chi(\epsilon x)\frac{f(v_n)}{v_n}w_n^2 + \big(1-\chi(\epsilon x)\big)\frac{\tilde{f}(v_n)}{v_n}w_n^2\right)dx,$$

and then

$$\limsup_{n\to\infty}\int_{\mathbb{R}^N}\chi(\epsilon x)\frac{f(v_n)}{v_n}w_n^2\,dx\leqslant 1,$$

which implies that

$$\limsup_{n\to\infty}\int_{\mathbb{R}^N}\chi(\epsilon x+\epsilon_n y_n)\frac{2F(v_n(x+y_n))}{v_n(x+y_n)^2}w_n(x+y_n)^2\,dx\leqslant 1. \tag{2.13}$$

By (2.12), $v_n(x+y_n)\to+\infty$ as $n\to\infty$, for all $x\in\Gamma$. Fatou's lemma implies that

$$\liminf_{n\to\infty}\int_{\Gamma}\chi(\epsilon x+\epsilon_n y_n)\frac{2F(v_n(x+y_n))}{v_n(x+y_n)^2}w_n(x+y_n)^2\,dx=\infty,$$

which contradicts (2.13). Therefore i) does not hold. On the other hand, assuming ii) note that $\chi(\epsilon_n\cdot)w_n\in H^1(\mathbb{R}^N)$ and $\{\chi(\epsilon_n\cdot)w_n\}$ is a bounded sequence in $H^1(\mathbb{R}^N)$. Then we can use Lemma I.1 [14] to conclude that

$$\|\chi(\epsilon_n\cdot)w_n\|_{L^r(\mathbb{R}^N)}\to 0,\quad \text{as } n\to\infty, \text{for all } 2<r<2_*. \tag{2.14}$$

Let $t_n\in[0,1]$ be such that

$$I_{\epsilon_n}(t_nv_n)=\max_{t\in[0,1]}I_{\epsilon_n}(tv_n).$$

As $\|v_n\|_{\epsilon_n}\to+\infty$, for a given $R>0$ there exists $n_0\in\mathbb{N}$ such that $R/\|v_n\|_{\epsilon_n}\leqslant 1$ for all $n\geqslant n_0$. Then, $(g_3)$ implies that

$$\begin{aligned}I_{\epsilon_n}(t_nv_n)&\geqslant I_{\epsilon_n}\left(\frac{R}{\|v_n\|_{\epsilon_n}}v_n\right)\\&\geqslant\frac{R^2}{2}\|w_n\|_{\epsilon_n}^2-\int_{\mathbb{R}^N}\chi(\epsilon_nx)F(Rw_n)\,dx-\frac{R^2}{2}\int_{\mathbb{R}^N}\frac{V(\epsilon_nx)}{k}w_n^2\,dx\\&\geqslant\frac{R^2(k-1)}{2k}-\int_{\mathbb{R}^N}\chi(\epsilon_nx)F(Rw_n)\,dx.\end{aligned}$$

By $(f_2)$ and $(f_3)$, for $\eta>0$,

$$\int_{\mathbb{R}^N}\chi(\epsilon_nx)F(Rw_n)\,dx\leqslant\eta R^2\|w_n\|_{L^2(\mathbb{R}^N)}^2+C_\eta R^{p+1}\int_{\mathbb{R}^N}\chi(\epsilon_nx)|w_n|^{p+1}\,dx.$$

From (2.14),

$$\begin{aligned}\int_{\mathbb{R}^N}\chi(\epsilon_nx)F(Rw_n)\,dx&\leqslant\eta R^2\|w_n\|_{E_\epsilon}^2+C_\eta R^{p+1}\|\chi(\epsilon_n\cdot)w_n\|_{L^{p+1}(\mathbb{R}^N)}\|w_n\|_{L^{p+1}(\mathbb{R}^N)}^p\\&\leqslant\eta R^2\|w_n\|_{E_\epsilon}^2+o_n(1)=\eta R^2+o_n(1).\end{aligned}$$

Since $\eta>0$ is arbitrary we have that

$$\liminf_{n\to+\infty}I_{\epsilon_n}(t_nv_n)\geqslant\frac{R^2(k-1)}{2k},$$

for all $R>0$ and then

$$\lim_{n\to+\infty}I_{\epsilon_n}(t_nv_n)=+\infty.$$

Since $I_{\epsilon_n}(0)=0$ and $I_{\epsilon_n}(v_n)\to d$, we have that $t_n\in(0,1)$. Then

$$I'_{\epsilon_n}(t_nv_n)v_n=0.$$

By $(g_4)$ we have that

$$\begin{aligned} 2I_{\epsilon_n}(t_n v_n) &= 2I_{\epsilon_n}(t_n v_n) - I'_{\epsilon_n}(t_n v_n)t_n v_n \\ &= \int_{\mathbb{R}^N} \big(g(\epsilon x, t_n v_n)t_n u_n - 2G(\epsilon x, t_n v_n)\big)\,dx \\ &\leqslant \int_{\mathbb{R}^N} \big(g(\epsilon x, v_n)u_n - 2G(\epsilon x, v_n)\big)\,dx, \end{aligned}$$

which implies that

$$\lim_{n\to+\infty} \int_{\mathbb{R}^N} \big(g(\epsilon x, v_n)v_n - 2G(\epsilon x, v_n)\big)\,dx = +\infty. \tag{2.15}$$

On the other hand

$$\int_{\mathbb{R}^N} \big(g(\epsilon_n x, v_n)v_n - 2G(\epsilon_n x, v_n)\big)\,dx = 2I_{\epsilon_n}(v_n) - I'_{\epsilon_n}(v_n)v_n = 2d + o_n(1),$$

which contradicts (2.15). Therefore the result follows. □

**Lemma 2.6.** *There exist $\{y_n\} \subset \mathbb{R}^N$ and $R, \beta > 0$ such that*

$$\liminf_{n\to\infty} \int_{B_R(y_n)} v_n^2\,dx \geqslant \beta > 0.$$

**Proof.** Suppose the assertion of the lemma is false. Then by Lemma I.1 of [14] (with $q = 2$ and $p = \frac{2N}{N-2}$), $v_n \to 0$ in $L^r(\mathbb{R}^N)$ where $2 \leqslant r \leqslant 2_*$. Hence by the Lebesgue Dominated Convergence Theorem, we get

$$\int_{\mathbb{R}^N} g(\epsilon_n x, v_n)v_n\,dx = o_n(1) \quad \text{and} \quad \int_{\mathbb{R}^N} G(\epsilon_n x, v_n)\,dx = o_n(1).$$

Then $c_{\epsilon_n} \to 0$ as $n \to \infty$. On the other hand, since the minimax value is an increasing function of the potential we have $c_{\epsilon_n} \geqslant d$, $\forall n \in \mathbb{N}$, where $d > 0$ is the minimax value associated to the problem

$$\Delta^2 v + V_0 v = g(\epsilon_n x, v) \quad \text{in } \mathbb{R}^N.$$

This contradiction proves the lemma. □

For $R > 0$ given by Lemma 2.6, we have:

**Lemma 2.7.** *The sequence $\{\epsilon_n y_{\epsilon_n}\}$ is bounded and $dist(\epsilon_n y_{\epsilon_n}, \Omega) \leqslant \epsilon_n R$.*

**Proof.** Let $K_\delta$ denote a $\delta$-neighborhood of $\Omega$, where $\delta > 0$. Let $\phi \in C^\infty(\mathbb{R}^N)$ be a cut-off function such that $\phi = 0$ in $\Omega$, $\phi = 1$ in $\mathbb{R}^N \backslash K_\delta$, $0 \leqslant \phi \leqslant 1$, $|\nabla\phi| \leqslant C/\delta$ and $|\Delta\phi| \leqslant C/\delta^2$. Setting $\phi_\epsilon(x) = \phi(\epsilon x)$ and using $v_{\epsilon_n}\phi_{\epsilon_n}$ as test function in (2.3) we have

$$\int_{\mathbb{R}^N} \big(\Delta v_{\epsilon_n} \Delta(v_{\epsilon_n}\phi_{\epsilon_n}) + V(\epsilon_n x)v_{\epsilon_n}^2 \phi_{\epsilon_n}\big)\,dx = \int_{\mathbb{R}^N} g(\epsilon_n x, v_{\epsilon_n})v_{\epsilon_n}\phi_{\epsilon_n}\,dx,$$

which gives

$$\begin{aligned} V_0\left(1 - \frac{1}{k}\right) \int_{\mathbb{R}^N} v_{\epsilon_n}^2 \phi_{\epsilon_n}\,dx &\leqslant \frac{C\epsilon_n}{\delta} \|v_{\epsilon_n}\|^2_{H^2(\mathbb{R}^N)} + \frac{C\epsilon_n^2}{\delta^2} \|v_{\epsilon_n}\|^2_{H^2(\mathbb{R}^N)} \\ &\leqslant \frac{C\epsilon_n}{\delta} \|v_{\epsilon_n}\|^2_{H^2(\mathbb{R}^N)}. \end{aligned}$$

If there is $\{\epsilon_k\}$ subsequence such that $B_R(y_{\epsilon_k}) \cap K_\delta/\epsilon_k = \emptyset$, then

$$\begin{aligned} V_0\left(1 - \frac{1}{k}\right) \int_{B_R(y_{\epsilon_k})} v_{\epsilon_k}^2\,dx &\leqslant V_0\left(1 - \frac{1}{k}\right) \int_{\mathbb{R}^N} v_{\epsilon_k}^2 \phi_{\epsilon_k}\,dx \\ &\leqslant \frac{C\epsilon_k}{\delta} \|v_{\epsilon_k}\|^2_{H^2(\mathbb{R}^N)} \to 0 \end{aligned}$$

as $k \to \infty$, which contradicts Lemma 2.6. Hence, for each $n \in \mathbb{N}$ there exists $x_n$ such that $\epsilon_n x_n \in K_\delta$ and $|y_{\epsilon_n} - x_n| < R$. Hence, $dist(\epsilon_n y_{\epsilon_n}, \Omega) < \epsilon_n R + \delta$, for all $\delta > 0$, which completes the proof. □

**Remark 2.1.** It is worth pointing out that by Lemma 2.7, we can assume that $\epsilon_n y_{\epsilon_n} \in \overline{\Omega}$ for all $n$ sufficiently large. In fact, on the contrary, we consider $\epsilon_n^{-1} z_n$ instead of $y_{\epsilon_n}$, where $z_n \in \overline{\Omega}$ is such that $|\epsilon_n y_n - z_n| < \epsilon_n R$. This fact will be used to guarantee that $\epsilon_n y_{\epsilon_n} \to x_0' \in \overline{\Omega}$.

The following result plays a central role in the proof of Theorem 1.1.

**Lemma 2.8.** *The following assertions hold*:

(i) $\lim_{n\to\infty} c_{\epsilon_n} = c_0$,
(ii) $\lim_{n\to\infty} V(\epsilon_n y_{\epsilon_n}) = V_0$.

**Proof.** By Lemma 2.7, we can assume that $\epsilon_n y_{\epsilon_n} \to x_0' \in \overline{\Omega}$ along a subsequence. Let us consider $w_n(x) = v_{\epsilon_n}(x + y_{\epsilon_n})$. Note that by Lemma 2.6

$$\liminf_{n\to\infty} \int_{B_R(0)} w_n^2 \, dx \geqslant \beta > 0.$$

Using that $\{v_{\epsilon_n}\}$ is bounded, it follows that there exists $w \in H^2(\mathbb{R}^N)\backslash\{0\}$, such that

$$\begin{aligned} &w_n \rightharpoonup w \quad \text{in } H^2(\mathbb{R}^N), \\ &w_n \to w \quad \text{in } L^r_{loc}(\mathbb{R}^N) \text{ where } 2 \leqslant r < 2_*, \\ &w_n \to w \quad \text{a.e. in } \mathbb{R}^N. \end{aligned}$$

Since $w_n$ satisfies

$$\Delta^2 w_n + V(\epsilon_n x + \epsilon_n y_n) w_n = g(\epsilon_n x + \epsilon_n y_n, w_n) \quad \text{in } \mathbb{R}^N, \tag{2.16}$$

we have that $w$ satisfies

$$\Delta^2 w + V(x_0') w = \chi(x) f(w) + (1 - \chi(x)) \tilde{f}(x) = \tilde{g}(x, w) \quad \text{in } \mathbb{R}^N, \tag{2.17}$$

where $\chi(x) = \lim_{n\to\infty} \chi_\Omega(\epsilon_n x + \epsilon_n y_n)$, almost everywhere in $\mathbb{R}^N$.

Denote by $\tilde{I}$ the energy functional associated with the problem (2.17) and by $\tilde{c}$ its minimax level. We now consider the problem

$$\Delta^2 w + V(x_0') w = f(w) \quad \text{in } \mathbb{R}^N \tag{2.18}$$

and let denote by $\bar{c}$ the minimax level of the functional $\bar{I}$ associated with the problem (2.18). In the following we show that $\bar{c} = \tilde{c}$. Since

$$\tilde{G}(x, s) = \int_0^s \tilde{g}(x, t) \, dt \leqslant F(s),$$

we have $\bar{I}(u) \leqslant \tilde{I}(u)$ for all $u \in H^2(\mathbb{R}^N)$, which implies that $\bar{c} \leqslant \tilde{c}$. In order to verify that $\tilde{c} \leqslant \bar{c}$, it is sufficient to prove that

$$\tilde{I}(w) \leqslant \liminf_{n\to\infty} c_{\epsilon_n}. \tag{2.19}$$

In fact, if (2.19) holds, then

$$\tilde{c} \leqslant \tilde{I}(w) \leqslant \liminf_{n\to\infty} c_{\epsilon_n} \leqslant \limsup_{n\to\infty} c_{\epsilon_n} \leqslant c_0 \leqslant \bar{c}. \tag{2.20}$$

Therefore, $\tilde{c} = \bar{c}$, which implies that $\tilde{c} = \bar{c} = c_0$. By (2.20), we have

$$\lim_{n\to\infty} c_{\epsilon_n} = c_0,$$

which proves (i).

The proof of (2.19) is based on some ideas of [10, Lemma 2.2]. By elliptic estimates we see that $w_n$ converges to $w$ in $C^2_{loc}(\mathbb{R}^N)$. Hence, for all $R > 0$, we have

$$\lim_{n\to\infty}\left[\int_{B_R(0)}\left(\frac{1}{2}\left(|\Delta w_n|^2+V(\epsilon_n x+\epsilon_n y_n)w_n^2\right)-G(\epsilon_n x+\epsilon_n y_n,w_n)\right)dx\right]$$
$$=\int_{B_R(0)}\left(\frac{1}{2}\left(|\Delta w|^2+V\left(x_0'\right)w^2\right)-\tilde{G}(x,w)\right)dx.$$

Since the last integral converges to $\tilde{I}(w)$ as $R\to+\infty$, for all $\delta>0$ there exists $R>0$ sufficiently large such that

$$\lim_{n\to\infty}\left[\int_{B_R(0)}\left(\frac{1}{2}\left(|\Delta w_n|^2+V(\epsilon_n x+\epsilon_n y_n)w_n^2\right)-G(\epsilon_n x+\epsilon_n y_n,w_n)\right)dx\right]$$
$$\geqslant\tilde{I}(w)-\delta.$$

Hence, in order to prove (2.19), it suffices to show that for $R>0$ sufficiently large, we have

$$\liminf_{n\to\infty}\left[\int_{B_R^c(0)}\left(\frac{1}{2}\left(|\Delta w_n|^2+V(\epsilon_n x+\epsilon_n y_n)w_n^2\right)-G(\epsilon_n x+\epsilon_n y_n,w_n)\right)dx\right]$$
$$\geqslant-\delta. \tag{2.21}$$

Consider a smooth cut-off function $\eta_R$ such that $\eta_R=1$ in $B_R^c(0)$, $\eta_R=0$ in $B_{R-1}(0)$ and $|\nabla\eta_R|,|\Delta\eta_R|\leqslant C$, where $C$ is independent of $R$. Using $\varphi=\eta_R w_n$ as a test function in (2.16), yields

$$0=\int_{\mathbb{R}^N}\left(\Delta w_n\Delta(\eta_R w_n)+V(\epsilon_n x+\epsilon_n y_n)\eta_R w_n^2-g(\epsilon_n x+\epsilon_n y_n,w_n)\eta_R w_n\right)dx$$
$$=A_{1,n}+A_{2,n}+A_{3,n},$$

where

$$A_{1,n}=\int_{B_R^c(0)}\left(|\Delta w_n|^2+V(\epsilon_n x+\epsilon_n y_n)w_n^2-2G(\epsilon_n x+\epsilon_n y_n,w_n)\right)dx,$$
$$A_{2,n}=\int_{B_R^c(0)}\left(2G(\epsilon_n x+\epsilon_n y_n,w_n)-g(\epsilon_n x+\epsilon_n y_n,w_n)w_n\right)dx$$

and

$$A_{3,n}=\int_{A_{R,R-1}}\left[\Delta w_n\Delta(\eta_R w_n)+V(\epsilon_n x+\epsilon_n y_n)(\eta_R w_n)^2-g(\epsilon_n x+\epsilon_n y_n,w_n)\eta_R w_n\right]dx$$

where $A_{R,R-1}$ is the annulus $B_R(0)\backslash B_{R-1}(0)$. By $(g_3)$, $A_{2,n}\leqslant 0$. On the other hand, we can choose $R>0$ sufficiently large such that

$$\lim_{n\to\infty}|A_{3,n}|\leqslant\delta.$$

Hence, as $A_{1,n}=-A_{2,n}-A_{3,n}$,

$$\liminf_{n\to\infty}\frac{1}{2}A_{1,n}\geqslant-\frac{\delta}{2}$$

which shows that (2.21) holds. Consequently (2.19) holds.

Finally, suppose by contradiction that (ii) is false. Thus, $V(x_0')>V_0$ and so $\tilde{c}>c_0$, which is impossible. This concludes the proof of Lemma 2.8. □

By the proof of Lemma 2.8, we see that $\epsilon_n y_n\to x_0'\in\overline{\Omega}$. Nevertheless, by the definition of $\Omega'$ and $(V_2)$, it follows that $x_0'\in\Omega'$. Indeed, if $x_0'\in\overline{\Omega}\backslash\Omega'$, since

$$\inf_{\overline{\Omega}\backslash\Omega'}V>\inf_{\mathbb{R}^N}V,$$

we have that $V(x_0')>V_0$, which contradicts Lemma 2.8ii).

## 3. Uniform decay

Although we have obtained existence and concentration of solutions to the modified problem, nothing can be said about the original one. In order to prove that the function $v_n$ is in fact a solution of the original problem, we will prove a kind of uniform decay at infinity of the translations $w_n(x) = v_{\epsilon_n}(x + y_{\epsilon_n})$. We first prove a technical lemma:

**Lemma 3.1.** *Let $\{u_n\} \subset H^2(\mathbb{R}^N)$ be a $(PS)_d$ sequence to $I_0$. If $u_n \rightharpoonup 0$ and $u_n \nrightarrow 0$ in $H^2(\mathbb{R}^N)$, then $d \geqslant c_0$.*

**Proof.** Let $s_n > 0$ be such that $s_n u_n \in \mathcal{N}_0$. Let us show that

$$\limsup_{n\to\infty} s_n \leqslant 1. \tag{3.1}$$

Suppose by contradiction that there exists $\delta > 0$ such that

$$s_n \geqslant 1 + \delta, \quad \forall n \in \mathbb{N}, \tag{3.2}$$

along a subsequence. Using the same arguments as that in Lemma 2.5, one can prove that $\{u_n\}$ is a bounded sequence in $H^2(\mathbb{R}^N)$. Since $I_0'(u_n)u_n = o_n(1)$ and $I_0'(s_n u_n)s_n u_n = 0$ it follows that

$$\int_{\mathbb{R}^N} \left( \frac{f(s_n u_n)}{s_n u_n} - \frac{f(u_n)}{u_n} \right) u_n^2 \, dx = o_n(1). \tag{3.3}$$

By Lemma I.1 in [14], there exist $R, \beta > 0$ and $\{y_n\} \subset \mathbb{R}^N$ such that

$$\liminf_{n\to\infty} \int_{B_R(y_n)} u_n^2 \, dx \geqslant \beta. \tag{3.4}$$

Let us set $\bar{u}_n(x) = u_n(x + y_n)$ and note that $\{\bar{u}_n\}$ is also a bounded sequence in $H^2(\mathbb{R}^N)$. Then, $\bar{u}_n \rightharpoonup \bar{u}$ in $H^2(\mathbb{R}^N)$ and by (3.4), $\bar{u} \neq 0$ in $\Lambda \subset B_{R_1}(0)$, with $|\Lambda| > 0$. By $(f_5)$, (3.2), (3.3) and Fatou's lemma, we have that

$$0 < \int_{\Lambda} \left( \frac{f((1+\delta)\bar{u})}{(1+\delta)\bar{u}} - \frac{f(\bar{u})}{\bar{u}} \right) \bar{u}^2 \, dx = 0,$$

which give us a contradiction. Then (3.1) holds.

Now we have two cases to consider:

- $s_n \to s < 1$ along a subsequence. Without loss of generality, let us suppose that $s_n < 1$ for all $n \in \mathbb{N}$. Then

$$\begin{aligned} c_0 &\leqslant I_0(s_n u_n) \\ &= \int_{\mathbb{R}^N} \left( \frac{1}{2} f(s_n u_n) s_n u_n - F(s_n u_n) \right) dx \\ &\leqslant \int_{\mathbb{R}^N} \left( \frac{1}{2} f(u_n) u_n - F(u_n) \right) dx \\ &= I_0(u_n) - \frac{1}{2} I_0'(u_n) u_n + o_n(1) \\ &= d + o_n(1) \end{aligned}$$

  and in this case we have the result.
- There exists a subsequence of $\{s_n\}$ such that $s_n \to 1$, as $n \to \infty$. In this case,

$$d + o_n(1) = I_0(u_n) = I_0(s_n u_n) + I_0(u_n) - I_0(s_n u_n).$$

  Then

$$d + o_n(1) \geqslant c_0 + I_0(u_n) - I_0(s_n u_n). \tag{3.5}$$

  Note that

$$I_0(u_n) - I_0(s_n u_n) = \frac{(1 - s_n^2)}{2} \int_{\mathbb{R}^N} \left( |\Delta u_n|^2 dx + V_0 u_n^2 \right) dx + \int_{\mathbb{R}^N} \left( F(s_n u_n) - F(u_n) \right) dx = o_n(1), \tag{3.6}$$

where the last equality can be verified by using the Mean Value Theorem. Therefore, (3.5) becomes

$$d + o_n(1) \geqslant c_0 + o_n(1),$$

which gives us the result. □

Now let us prove a compactness result in Nehari manifolds that will be used to prove the uniform decay.

**Lemma 3.2.** *Let $\{z_n\} \subset H^2(\mathbb{R}^N)$ be such that $I_0(z_n) \to c_0$ and $z_n \in \mathcal{N}_0$, for all $n \in \mathbb{N}$. If $z_n \rightharpoonup z \neq 0$, then $z_n \to z$ in $H^2(\mathbb{R}^N)$ along a subsequence.*

**Proof.** By the Ekeland Variational Principle, we can assume that $\{z_n\}$ is a $(PS)_{c_0}$ sequence for $I_0$ in $H^2(\mathbb{R}^N)$. Then it is possible to show that $I_0'(z) = 0$, which implies that $z \in \mathcal{N}_0$. By Fatou's lemma,

$$\begin{aligned} c_0 &= \lim_{n\to\infty}\left[\int_{\mathbb{R}^N}\left(\frac{1}{2}f(z_n)z_n - F(z_n)\right)dx\right] \\ &\geqslant \int_{\mathbb{R}^N}\left(\frac{1}{2}f(z)z - F(z)\right)dx = I_0(z) \geqslant c_0, \end{aligned}$$

which implies that

$$I_0(z) = c_0. \tag{3.7}$$

Let $u_n = z_n - z$ and note that by Brezis–Lieb lemma, $\{u_n\}$ is $(PS)_d$ sequence for $I_0$ where $d = c_0 - I_0(z) = 0$. Note that $u_n \rightharpoonup 0$ in $H^2(\mathbb{R}^N)$ and we claim that in fact $u_n \to 0$ in $H^2(\mathbb{R}^N)$. Indeed, suppose by contradiction that $u_n \nrightarrow 0$. By Lemma 3.1 we have that $d \geqslant c_0 > 0$, which give us a contradiction with the fact that $d = 0$. Therefore we have the result. □

**Lemma 3.3.** *The sequence $\{w_n\}$ contains a strongly convergent subsequence in $H^2(\mathbb{R}^N)$.*

**Proof.** By Lemma 2.8, we have

$$\lim_{n\to\infty} I_{\epsilon_n}(v_n) = \lim_{n\to\infty} c_{\epsilon_n} = c_0.$$

Denote by $\mathcal{N}_0$ the Nehari manifold associated to (2.6). Given $v \in H^2(\mathbb{R}^N)\backslash\{0\}$, from $(g_4)$, there exists $\varphi_0(v) > 0$ such that $\varphi_0(v)v \in \mathcal{N}_0$. Set $\tilde{w}_n = \varphi_0(w_n)w_n$. Hence,

$$\begin{aligned} c_0 &\leqslant \frac{1}{2}\int_{\mathbb{R}^N}\left(|\Delta\tilde{w}_n|^2 + V_0\tilde{w}_n^2\right)dx - \int_{\mathbb{R}^N}F(\tilde{w}_n)\,dx \\ &\leqslant \frac{1}{2}\int_{\mathbb{R}^N}\left(|\Delta\tilde{w}_n|^2 + V(\epsilon_n x + \epsilon_n y_n)\tilde{w}_n^2\right)dx - \int_{\mathbb{R}^N}G(\epsilon_n x + \epsilon_n y_n, \tilde{w}_n)\,dx \\ &= I_{\epsilon_n}\big(\varphi_0(w_n)v_n\big) \leqslant I_{\epsilon_n}(v_n) = c_{\epsilon_n} = c_0 + o_n(1), \end{aligned}$$

which implies that $I_0(\tilde{w}_n) \to c_0$ as $n \to \infty$.

We now prove that $\varphi_0(w_n) \to \varphi_0 > 0$ along a subsequence. We first observe that there exists $M > 0$ such that $|\varphi_0(w_n)| \leqslant M$, $\forall n \in \mathbb{N}$. In fact, since $w_n \nrightarrow 0$ there exists $\delta > 0$ such that $\|w_n\|_{H^2(\mathbb{R}^N)} > \delta$ along a subsequence. On the other hand, it is easy to see that $\{\tilde{w}_n\}$ is a bounded sequence in $H^2(\mathbb{R}^N)$. Then

$$\big|\varphi_0(w_n)\big|\delta < \big\|\varphi_0(w_n)w_n\big\|_{H^2(\mathbb{R}^N)} \leqslant K$$

which implies that

$$\big|\varphi_0(w_n)\big| \leqslant \frac{K}{\delta} = M, \quad \forall n \in \mathbb{N}.$$

Hence, $\varphi_0(w_n) \to \varphi_0 \geqslant 0$. We now observe that $\varphi_0 > 0$, otherwise

$$\|\tilde{w}_n\|_{H^2(\mathbb{R}^N)} = \big|\varphi_0(w_n)\big|\|w_n\|_{H^2(\mathbb{R}^N)} \to 0$$

as $n \to \infty$, which is impossible. Therefore $\tilde{w}_n = \varphi_0(w_n)w_n \rightharpoonup \varphi_0 w \neq 0$ in $H^2(\mathbb{R}^N)$. Therefore, we conclude from Lemma 3.2 that the lemma is proved. □

Combining Lemma 3.3 with the Sobolev imbeddings, it follows that $w_n \to w$ in $L^{2_*}(\mathbb{R}^N)$. Therefore, we obtain

$$\int_{B_R^c(0)} |w_n|^{2_*}\,dx \to 0 \quad \text{as } R \to \infty \text{ uniformly in } n. \tag{3.8}$$

**Lemma 3.4.** *$w_n(x) \to 0$ as $|x| \to \infty$, uniformly in n.*

**Proof.** By the uniform $L^\infty$ estimates to solutions of subcritical biharmonic equations given in [15], we have

$$\|w_n\|_{L^\infty(\mathbb{R}^N)} \leqslant C, \quad \forall n \in \mathbb{N},$$

where $C$ is independent of $n$. Given any $x \in \mathbb{R}^N$, the function $w_n \in L^q(B_1(x))$ for all $q \geqslant 1$. By [1, Theorem 7.1] it follows that

$$\begin{aligned}\|w_n\|_{W^{4,q}(B_1(x))} &\leqslant C\big(\big\|f(w_n)\big\|_{L^q(B_2(x))} + \|w_n\|_{L^q(B_2(x))}\big) \leqslant C\|w_k\|_{L^q(B_2(x))}\\ &\leqslant C\|w_k\|_{L^\infty(\mathbb{R}^N)}^{\frac{q-2_*}{q}}\|w_k\|_{L^{2_*}(B_2(x))}^{2_*} = C\|w_k\|_{L^{2_*}(B_2(x))}^{2_*},\end{aligned}$$

with $C > 0$ being a constant independent of $x$ and $n$. If $q > N$, we have the continuous imbedding $W^{4,q}(B_1(x)) \hookrightarrow C^{3,\alpha}(\overline{B_1(x)})$ for $\alpha \in (0, 1 - \frac{N}{q})$. Then

$$\|w_k\|_{C^{3,\alpha}(\overline{B_1(x)})} \leqslant \|w_k\|_{W^{4,q}(B_1(x))} \leqslant C\|w_k\|_{L^{2_*}(B_2(x))}^{2_*}.$$

By (3.8), it follows that $|w_n(x)| \to 0$ as $|x| \to \infty$ uniformly in $n$. □

Finally, we are ready to prove that $v_n$ is in fact a solution of (2.1). Let $n_0 \in \mathbb{N}$ and $\rho > 0$ be such that

$$\big|w_n(x)\big| < a, \quad \forall x \in B_\rho(0)^c,\ \forall n \geqslant n_0.$$

Since $x_0' \in \Omega'$ and $\epsilon_n y_n \to x_0'$, it is possible to choose $n_1 \in \mathbb{N}$ such that $B_\rho(0) \subset (\Omega'_{\epsilon_n} - y_n)$, for all $n \geqslant n_1$. Taking $n \geqslant \max\{n_0, n_1\}$, we have

$$g\big(\epsilon_n x + \epsilon_n y_n, w_n(x)\big) = f\big(w_n(x)\big), \quad \forall x \in \mathbb{R}^N.$$

Hence, for $n \geqslant \max\{n_0, n_1\}$ it follows that $w_n$ satisfies

$$\Delta^2 w_n + V(\epsilon_n x + \epsilon_n y_n) w_n = f(w_n) \quad \text{in } \mathbb{R}^N,$$

which implies that $v_n$ satisfies (2.1).

In order to prove the concentration behavior of solutions, we claim that there exists $\rho > 0$ such that $\|u_n\|_{L^\infty(\mathbb{R}^N)} = \|w_n\|_{L^\infty(\mathbb{R}^N)} > \rho$, for all $n \in \mathbb{N}$ along a subsequence. In fact, if $\|w_n\|_{L^\infty(\mathbb{R}^N)} \to 0$, then

$$\begin{aligned}\|w_n\|_{H^2(\mathbb{R}^N)}^2 &\leqslant C\int_{\mathbb{R}^N}\big(|\Delta w_n|^2 + V(\epsilon_n x + \epsilon_n y_n)w_n^2\big)\,dx\\ &= C\int_{\mathbb{R}^N} f(w_k)w_k\,dx\\ &\leqslant \|w_k\|_{L^\infty(\mathbb{R}^N)}\int_{\mathbb{R}^N} f(w_k)\,dx \to 0\end{aligned}$$

as $n \to \infty$, which contradicts the fact that $w_n \to w$ and $w \neq 0$.

Let $x_n$ be the maximum point of $|u_n|$ in $\mathbb{R}^N$, then

$$p_n := \frac{x_n - \epsilon_n y_n}{\epsilon_n}$$

is the maximum point of $|w_n|$. By Lemma 3.4, there exists $R_0 > 0$ such that $p_n \in B_{R_0}(0)$ for all $n$ sufficiently large. Then, along a subsequence $p_n \to p_0$ as $n \to \infty$. Hence

$$x_n = \epsilon_n p_n + \epsilon_n y_n \to x_0' \in \Omega', \quad \text{as } n \to \infty,$$

where $V(x_0') = V_0$, which proves Theorem 1.1.

## Acknowledgments

The authors are grateful to Claudianor O. Alves and Marco A. S. Souto for valuable discussions.